\documentclass[11pt]{article}

\usepackage{latexsym,amssymb}
\usepackage{enumerate}
\usepackage{amsmath}
\usepackage[latin1]{inputenc}
\usepackage[english]{babel}

\newtheorem{theorem}{Theorem}[section]

\newtheorem{lemma}[theorem]{Lemma}

\newenvironment{proof}{\trivlist\item[\hskip
            \labelsep{\it
Proof:\/}]\ignorespaces}{\hfill$\Box$\endtrivlist}
\newcommand{\bg}{\begin}
\newcommand{\e}{\end}
\newcommand{\ee}{\e{enumerate}}
\newcommand{\be}{\bg{enumerate}}
\newcommand{\ei}{\e{itemize}}
\newcommand{\bi}{\bg{itemize}}
\newcommand{\MV}[1]{\langle #1, \oplus, \lnot,0
\rangle }

\begin{document}
\title{ A note on bases of admissible rules of proper axiomatic extensions of $\L$ukasiewicz logic}

\author{\textsc{Joan Gispert}\\
Universitat de Barcelona}
\date{\texttt{email:jgispertb@ub.edu}}

\maketitle

\begin{abstract}
\noindent In this note we prove that single-conclusion admissible rules of any  proper axiomatic extension of the infinite valued $\L$ukasiewicz logic are finitely based. The proof strongly relies on the characterization of least $\mathbf{V}$-quasivarieties given in \cite{gispert14}.
\end{abstract}

\section*{ Introduction.}

Admissible rules of a logic are those rules under which the set of theorems are closed. If  $L$ is a logic, an \emph{$L$-unifier} of a formula $\varphi$ is a substitution $\sigma$ such that $\vdash_{L}\sigma\varphi$. A \emph{single-conclusion rule} is an expression of the form $\Gamma/\varphi$ where $\varphi$ is a formula and $\Gamma$ is a finite set of formulas. As usual $\Gamma/\varphi$ is derivable in $L$ iff $\Gamma\vdash_{L}\varphi$. The rule $\Gamma/\varphi$ is \emph{admissible} in $L$ iff every common $L$-unifier of $\Gamma$ is also an $L$-unifier of $\varphi$. $\Gamma/\varphi$ is \emph{passive $L$-admissible}   iff $\Gamma$ has no common $L$-unifier. We say that a logic is \emph{structurally complete} iff every admissible rule is a derivable rule. Roughly speaking, a logic is structurally complete
iff every proper finitary extension must contain new theorems, as opposed
to nothing but new rules of inference (see for instance \cite{berg, lor, ryb}). Every logic $L$ has a unique structurally complete extension $L'$ with same theorems of $L$ \cite{berg}. In particular, structural completeness can be seen as a kind of maximality condition on a logic. We say that a logic is \emph{almost structurally complete} iff every admissible rule is either derivable rule or passive. It follows from a result of Dzik \cite{Dzi} that every finite valued $\L$ukasiewicz logic is almost structurally complete. Je$\check{\mbox{r}}$\'{a}bek uses this result \cite[Corollary 3.7]{Jer10a} to obtain for every $n>2$, one  rule that axiomatizes admissible rules of the $n$-valued $\L$ukasiewicz logic. In \cite{ Jer10b} the same  author proves that admissible rules of the infinite valued $\L$ukasiewicz calculus $\L_{\infty}$ are not finitely based, moreover he explicitly constructs an infinite  base of single-conclusion admissible rules for $\L_{\infty}$. The purpose of this work is to obtain a basis of single-conclusion rules of every proper axiomatic extension of $\L_{\infty}$. Admissibility theory  normally uses  proof-theoretic techniques, however in this case we will take an algebraic approach taking advantage of the algebraization of $\L_{\infty}$. In fact for algebraizable logics there is a analogous algebraic notion of admissible quasiequations and structurally complete and almost structurally complete quasivarieties (see for instance \cite{berg, ryb}).

 It is well known that $\L_{\infty}$  is algebraizable  and the class of MV-algebras $\mathbf{MV}$ is its equivalent  algebraic quasivariety semantics \cite{RTV}. It follows from the algebraization, that quasivarieties of $\mathbf{MV}$ are in 1-1 correspondence with finitary extensions of $\L_{\infty}$. Actually, there is a dual isomorphism from the lattice of all quasivarieties of $\mathbf{MV}$ and the lattice of all finitary extensions of $\L_{\infty}$. Moreover if we restrict this correspondence to varieties of $\mathbf{MV}$ we get the dual isomorphism from the lattice of all varieties of $\mathbf{MV}$ and the lattice of all axiomatic extensions of $\L_{\infty}$. Given an axiomatic extension $L$ of $\L_{\infty}$ whose equivalent variety semantics is $\mathbf{V}_{L}$, single-conclusion admissible rules of $L$ can be seen as valid quasiequations in the $\mathbf{V}_{L}$-free algebra under a countable set of generators $F_{\mathbf{V}_{L}}(\omega)$ in the following sense: $\gamma_{1}, \ldots, \gamma_{n}/\varphi$ is $L$-admissible if and only if $\gamma_{1}\approx 1\&\cdots \& \gamma_{n}\approx 1 \Rightarrow \varphi \approx 1$ is valid in $F_{\mathbf{V}_{L}}(\omega)$ \cite{ryb}. Hence the algebraic study of admissible rules of $L$ is the study of $\mathcal{Q}(F_{\mathbf{V}_{L}}(\omega))$ the quasivariety generated by $F_{\mathbf{V}_{L}}(\omega)$. The quasivariety $\mathcal{Q}(F_{\mathbf{V}_{L}}(\omega))$ is the least quasivariety that generates $\mathbf{V}_{L}$ as a variety. These quasivarieties were studied in \cite{gispert14}, where a Komori's type characterization is accomplished. In this paper we will prove that $\mathcal{Q}(F_{\mathbf{V}_{L}}(\omega))$ is the class of bipartite algebras of the least quasivariety generated by MV-chains that generates $\mathbf{V}_{L}$ as a variety. We use this result to prove that admissible rules for any proper axiomatic extension of $\L_{\infty}$ are finitely based and finally we give an effective axiomatization for each one.


\noindent The paper is organized as follows. First, in Section 1 we introduce the necessary definitions, notation and preliminary results on Universal Algebra and MV-algebra Theory that we use throughout the paper. Section 2 is first devoted to survey already existing results on varieties and quasivarieties of MV-algebras and later to obtain principal algebraic results:   Theorems \ref{BipQ} and \ref{AllStComp}. Finally, Section 3 contains  the main theorem,  Theorem \ref{Ax}, where we built a finite basis for each axiomatic extension of $\L_{\infty}$.

\section{ Definitions and first properties.}

We assume the reader is familiar with Universal Algebra \cite{bursan81, gratzer}. To fix notation
we denote by  $ \mathcal{I}$, $\mathcal{H}$, $\mathcal{S}$, $\mathcal{P}$,
$\mathcal{P}_{R}$ and $\mathcal{P}_{U}$ the   operators {\em isomorphic
image, homomorphic image, substructure, direct product,
reduced product and ultraproduct} res\-pectively. We recall
that a class $\mathbf{K}$ of algebras is a {\em variety} if and only if it is
closed under $\mathcal{H}$, $\mathcal{S}$ and $\mathcal{P}$. A class
$\mathbf{K}$ of algebras is a {\em quasivariety} if and only if it is
closed under $\mathcal{I}$, $\mathcal{S}$ and $\mathcal{P}_{R}$, or equivalently, under
$\mathcal{I}$, $\mathcal{S}$, $\mathcal{P}$ and $\mathcal{P}_{U}$.
A class
$\mathbf{K}$ of algebras is a {\em universal class} if and only if it is
closed under $\mathcal{I}$, $\mathcal{S}$ and $\mathcal{P}_{U}$.
Given a class $\mathbf{K}$ of algebras, the \emph{variety generated by} $\mathbf{K}$, denoted by $\mathcal{V}(\mathbf{K})$, is the least variety containing
$\mathbf{K}$. Similarly, the \emph{quasivariety generated by a class}
$\mathbf{K}$, which we denote by $\mathcal{Q}(\mathbf{K})$, is the least quasivariety
containing $\mathbf{K}$. $\mathcal{U}(\mathbf{K})$ denotes the \emph{universal class generated by }$\mathbf{K}$, that is  the least universal class containing $\mathbf{K}$.
We also recall that a class $\mathbf{K}$ of algebras is a
variety if and only if it is an equational class; $\mathbf{K}$ is a
quasivariety if and only if it is a quasiequational class; $\mathbf{K}$ is a universal class if and only if $\mathbf{K}$ is definable by first order universal sentences.

%

\bigskip

\noindent An {\em MV-algebra} is an algebra $\langle A,\oplus,\lnot,0
\rangle$ satisfying the following equations:
\begin{description}
\item{MV1} $\;\;(x\oplus y)\oplus z\approx x\oplus (y\oplus z)$
\item{MV2} $\;\;x\oplus y \approx y\oplus x$
\item{MV3} $\;\;x\oplus 0 \approx x$
\item{MV4} $\;\;\lnot(\lnot x)\approx x$
\item{MV5} $\;\;x\oplus \lnot 0\approx\lnot 0$
\item{MV6} $\;\;\lnot(\lnot x \oplus y) \oplus y \approx \lnot(\lnot y
\oplus x) \oplus x$.
\end{description}

We write $1$ instead of $\neg 0,$
$x \odot y$ instead of $\neg(\neg x \oplus \neg y)$   and $a\to b$ instead of $\lnot a \oplus b$.
Further,  for all $n\in\omega$, where $\omega$ is the set of all natural numbers,  and $x \in A$,
the  MV-operations $nx$ and $x^{n}$
are inductively defined by

$$0x=0, \,\,\,\,\, (n+1)x=x\oplus (nx)$$
and

$$x^{0}=1, \,\,\,\,\, x^{n+1}=x\odot(x^{n}).$$

\noindent Following tradition we assume that the
operation  $x^{n}$
takes precedence over  any other operation;
also  $\neg$  takes precedence over $\odot$,
$\odot$  takes precedence over $\oplus$, and
$\oplus$ takes precedence over $\to$.

\noindent As shown by Chang \cite{cha58}, for any MV-algebra $A$,
the stipulation
$a\le b$ iff
$a\to b= 1$ endows $A$ with a bounded distributive lattice-order $\langle
A,\lor, \land, 0, 1
\rangle$, called the {\em natural order} of $A$.

$$\;x \lor y =_{def} \;\lnot(\lnot x \oplus y) \oplus y.$$

$$\;x \land y =_{def} \;\lnot(\lnot x \lor\lnot y).$$

\noindent An MV-algebra whose natural order is total is said to be an {\em MV-chain}.

\medskip

\noindent We recall that a {\em lattice-ordered abelian group} (for short,
$\ell$-group) is an algebra $\langle G,\land ,\lor ,+,-,0\rangle$
such that $\langle
G,\land ,\lor\rangle$ is a lattice,
$\langle G,+,-,0\rangle$ is an abelian group and satisfies the following
equation:
 $$(x\lor y)+z \approx (x+z)\lor(y+z)
$$

\noindent For any  $\ell$-group $G$ and  element $0 < u \in G$, let
$\Gamma(G,u) = \MV{[0, u]}$ be defined by
$$
\label{equation:gamma} [0, u] = \{a\in G \,|\, 0\leq a\leq u\}, \;\;
a\oplus b = u\land(a+b),
\;\; \neg a=u-a.$$
 Then,
$\MV{[0, u]}$ is an MV-algebra. Further, for any $\ell$-groups
$G$ and $H$ with elements $0 < u \in G$ and
$0 < v \in H,$ and any
$\ell$-group homomorphism
$f: G \to H$ such that $f(u) = v$, let
$\Gamma(f)$ be the restriction of $f$ to $[0,u]$.
 An element $0 < u \in G$ is called a {\it strong unit}
iff for each $x
\in G$ there is an integer
$n \geq 1$ such that
$x \leq nu$. Then, as proved in \cite{mun86}, (see also \cite{CDM})
$\Gamma $ is a categorical equivalence from the category of
$\ell$-groups with strong unit, with
$\ell$-homomorphisms that preserve strong units, onto the category of
MV-algebras with
MV-homomorphisms.
Moreover the functor $\Gamma $ preserves embeddings and epimorphisms.

\bigskip

\noindent The following MV-algebras play an important role in the paper.
\bg{itemize}
\item $\displaystyle\mathrm{[0,1]}=\Gamma(\mathbb{R},1),$ where $\mathbb{R}$ is the totally
ordered group of the reals.
\item $\displaystyle \mathrm{[0,1]}\cap \mathbb{Q}=\Gamma(\mathbb{Q},1)=\MV{\{\frac{k}{m} : k\leq m<\omega\}}$,
  where $\mathbb{Q}$ is the totally ordered abelian group of the
  rationals.
\ei
For every $0<n<\omega $
\bi
\item $\displaystyle  \mathrm{L}_{n}=\Gamma(\mathbb{Z},n)=\MV{\{0, 1,\ldots, n\}}$,
 where $\mathbb{Z}$ is the totally ordered group of all integers. Notice that  $\mathrm{L}_{n}$ is isomorphic to the subalgebra
of  $\mathrm{[0,1]}$ given by $\{0, \frac{1}{n},  \frac{2}{n}, \ldots, \frac{n-1}{n}, 1\}$.
\item  $\displaystyle \mathrm{L}_{n}^{\omega}=\Gamma(\mathbb{Z}\times_{lex}
 \mathbb{Z},(n,0))=\MV{\{(k,i) :(0,0)\leq(k,i)\leq (n,0)\}}$,
where $\times_{lex}$ denotes the lexicographic product.
\item  $\displaystyle \mathrm{L}_{n}^{s}=\Gamma(\mathbb{Z}\times_{lex}
 \mathbb{Z},(n,s))=\MV{\{(k,i) :(0,0)\leq(k,i)\leq (n,s)\}}$, where $s\in \mathbb{Z}$ such that $0\leq s <n$. Notice that $\mathrm{L}_{n}^{\omega}=\mathrm{L}_{n}^{0}$.


\e{itemize}

\bigskip

\noindent As usual, by an {\it ideal} of an MV-algebra $A$
we mean the kernel $I$ of a homomorphism
$h$ of $A$ into some MV-algebra $B$.
In other words,
$0\in I$, $I$ is closed under
the $\oplus$ operation, and $x \leq y \in I$
implies $x \in I$. We denote by $I(A)$ the set of all ideals of $A$.

%
An ideal is {\it prime} iff it is
the kernel of a homomorphism of
$A$ into an MV-chain. We denote by $Spec(A)=\{I\in I(A) : I \mbox{ is prime}\}$.
An ideal is {\it maximal} iff
it is the kernel of a homomorphism
of $A$ into $\mathrm{[0,1]}$. We denote by $\mathcal{M}(A)=\{I\in I(A): I \mbox{ is maximal}\}$.  The \emph{radical of }$A$ denoted by $Rad(A)$ is the intersection of all maximal ideals of $A$. Notice that when $A$ is an MV-chain $\mathcal{M}(A)=\{Rad(A)\}$ and $Rad(A)=\{a\in A : a^{k}\neq 0 \mbox{ for all } k>0\}$.

An MV-algebra is said to be \emph{bipartite} iff there is $I\in I(A)$ such that $A/I\cong\mathrm{L}_{1}$.

\section{Varieties and quasivarieties of MV-algebras}

Since the class of all MV-algebras is definable by a set of equations, it is a variety that we denote by $\mathbf{MV} $. By Chang's Completeness
Theorem \cite{cha59} (see also \cite{CDM}),  $\mathbf{MV}$ is the variety  generated by the MV-algebra $\mathrm{[0,1]}$ (or $\mathrm{[0,1]}\cap \mathbb{Q}$), in symbols,
$$\mathbf{MV} = \mathcal{V}(\mathrm{[0,1]})=\mathcal{V}(\mathrm{[0,1]}\cap \mathbb{Q}).$$

\noindent Proper subvarieties of $\mathbf{MV}$ are well known. Komori proves the following characterization


\bg{theorem}{\rm \cite[Theorem 4.11]{kom81}}
\label{varietats} $\mathbf{V}$ is a proper subvariety of $\mathbf{MV}$ if and only if
there exist two disjoint finite subsets $I,J$ of positive integers, not both empty such that
\[\qquad \qquad  \mathbf{V}=\mathcal{V}(\{\mathrm{L}_{i} : i\in I\}\cup\{\mathrm{L}_{j}^{\omega} : j\in J\}).\qquad \qquad \Box\]
\e{theorem}

\noindent A pair $(I,J)$ of finite subsets of positive integers, not both empty is said to be \emph{reduced } iff for every $n\in I$, there is no $k\in (I\smallsetminus\{n\})\cup J$ such that $n|k$ and for every $m\in J$, there is no $k\in J\smallsetminus\{m\}$ such that $m|k$. In \cite{Pa99} Panti shows that there is a 1-1 correspondence between proper subvarieties of $\mathbf{MV}$ and reduced pairs of finite subsets of positive integers not both empty. Given a reduced pair $(I,J)$, we denote by $\mathcal{V}_{I,J}$ its associated subvariety. Moreover for every reduced pair $(I,J)$ there is a single equation in just one variable of the form $\alpha_{I,J}(x)\approx 1$ axiomatizing $\mathcal{V}_{I,J}$.

Quasivarieties of MV-algebras have been studied by this author in \cite{gismun, gismuntor, gispert, gispert14}. Particularly in \cite{gispert}, the author finds a characterization, classification an axiomatization of every quasivariety generated by MV-chains. Moreover  he obtains  necessary condition for finitely axiomatization that yields to the following result:

\begin{theorem}\label{PropQuasFinAx}
Every quasivariety generated by MV-chains contained in a proper subvariety of $\mathbf{MV}$ is finitely axiomatizable
\end{theorem}

Let $\mathbf{V}$ a variety of any type of algebras, a quasivariety $\mathbf{K}$ of same type  is a \emph{$\mathbf{V}$-quasivariety} provided that $\mathcal{V}(\mathbf{K})=\mathbf{V}$.  It follows from the work in  \cite{gispert} that $\mathcal{Q}_{I,J}^{1}:=\mathcal{Q}( \{\mathrm{L}_{m} \mid m\in I\}\cup\{\mathrm{L}_{n}^{1}\mid n\in J\})$ is the least $\mathcal{V}_{I,J}$-quasivariety generated by chains. However $\mathcal{Q}_{I,J}^{1}$ is not the least $\mathcal{V}_{I,J}$-quasivariety. In fact, for any  variety $\mathbf{V}$, not necessarily of MV-algebras, the  least $\mathbf{V}$-quasivariety is  $\mathcal{Q}(F_{\mathbf{V}}(\omega))$. In \cite{gispert14} we study least $\mathbf{V}$-quasivarieties and we obtain the following characterization.

\begin{theorem}{\rm\cite[Theorem 4.8]{gispert14}}
Let $(I,J)$ be a reduced pair. Then
$\mathcal{Q}_{I,J}:=\mathcal{Q}( \{\mathrm{L}_{1}\times\mathrm{L}_{m} \mid m\in
I\}\cup\{\mathrm{L}_{1}\times\mathrm{L}_{n}^{1}\mid n\in J\})$ is the least $\mathcal{V}_{I,J}$-quasivariety.
\end{theorem}

Next result establishes the relation between least $\mathbf{V}$-quasivarieties of MV-algebras and least $\mathbf{V}$-quasivarieties generated by chains.

\begin{theorem} \label{BipQ}
$\mathcal{Q}_{I,J}$ is the class of all bipartite algebras in $\mathcal{Q}_{I,J}^{1}$
\end{theorem}
\proof Let $BPQ^{1}_{I,J}$ be the class of all all bipartite algebras in $\mathcal{Q}_{I,J}^{1}=\mathcal{Q}( \{\mathrm{L}_{m} \mid m\in I\}\cup\{\mathrm{L}_{n}^{1}\mid n\in J\})$. Since being bipartite is preserved under $\mathcal{I}$, $\mathcal{S}$, $\mathcal{P}$ and $\mathcal{P}_{U}$ and  $\{\mathrm{L}_{1}\times\mathrm{L}_{m} \mid m\in
I\}\cup\{\mathrm{L}_{1}\times\mathrm{L}_{n}^{1}\mid n\in J\}\subseteqq BPQ^{1}_{I,J}$, then $\mathcal{Q}_{I,J}=\mathcal{I}\mathcal{S}\mathcal{P}\mathcal{P}_{U}(\{\mathrm{L}_{1}\times\mathrm{L}_{m} \mid m\in
I\}\cup\{\mathrm{L}_{1}\times\mathrm{L}_{n}^{1}\mid n\in J\})\subseteqq BPQ^{1}_{I,J}$.

\noindent In order to prove the other implication, let $A\in BPQ^{1}_{I,J}$. Let $\Delta=\{I\in Spec(A) : A/I\in \mathcal{Q}( \{\mathrm{L}_{m} \mid m\in I\}\cup\{\mathrm{L}_{n}^{1}\mid n\in J\})\}$. Since $\mathcal{Q}( \{\mathrm{L}_{m} \mid m\in I\}\cup\{\mathrm{L}_{n}^{1}\mid n\in J\})$ is relative congruence distributive quasivariety, $\Delta$ induces a natural subdirect representation of $A$: $A\hookrightarrow_{SD}\prod_{I\in\Delta}A/I : a\mapsto(a/I)_{I\in\Delta}$ where each $A/I$ is an MV-chain of $\mathcal{Q}( \{\mathrm{L}_{m} \mid m\in I\}\cup\{\mathrm{L}_{n}^{1}\mid n\in J\})$. Since $A$ is bipartite there exists $I\in \Delta$ such that $A/I\cong \mathrm{L}_{1}$.  Thus $A\in \mathcal{I}\mathcal{S}\mathcal{P}(\{\mathrm{L}_{1}\times B \mid B\in K\})$ where $K=\{A/I \mid I\in \Delta \mbox{ and } A/I\not\cong \mathrm{L}_{1}\}$. Since  $B$ is an MV-chain of $\mathcal{Q}^{1}_{I,J}$, then $B\in \mathcal{I}\mathcal{S}\mathcal{P}_{U}(\{\mathrm{L}_{m} \mid m\in I\}\cup\{\mathrm{L}_{n}^{1}\mid n\in J\})=\bigcup_{m\in I}\mathcal{I}\mathcal{S}(\mathrm{L}_{m}) \cup \bigcup_{n\in J}\mathcal{I}\mathcal{S}\mathcal{P}_{U}(\mathrm{L}_{n}^{1})$.
Thus for every $B\in K$, $\mathrm{L}_{1}\times B\in \bigcup_{m\in I}\mathcal{I}\mathcal{S}(\mathrm{L}_{1}\times\mathrm{L}_{m}) \cup \bigcup_{n\in J}\mathcal{I}\mathcal{S}\mathcal{P}_{U}(\mathrm{L}_{1}\times\mathrm{L}_{n}^{1})\subseteq \mathcal{Q}_{I,J}$. \hfill $\Box$

\bigskip

In \cite{MR}, G. Metcalfe  and C. R\"{o}thlisberger give the following characterization of almost structurally complete quasivarieties
\begin{theorem}{\rm \cite[Theorem 4.10]{MR}} Let $\mathbf{K}$ be a quasivariety.  The following are equivalent for any $B\in \mathcal{S}(F_{\mathbf{K}}(\omega))$
     \be
      \item $\mathbf{K}$ is almost structurally complete.
      \item $\mathcal{Q}(\{A\times B : A\in \mathbf{K}\})=\mathcal{Q}(F_{\mathbf{K}}(\omega))$.
      \item $\{A\times B : A\in \mathbf{K}\}\subseteqq\mathcal{Q}(F_{\mathbf{K}}(\omega))$.
     \ee
\end{theorem}


   \begin{theorem} \label{AllStComp} Let $(I,J)$ be a reduced pair. Then $\mathcal{Q}(\{\mathrm{L}_{n} : n\in I\}\cup\{\mathrm{L}_{m}^{1} : m\in J\})$ is almost structurally complete.
   \end{theorem}

   \proof Since $F_{\mathcal{V}_{I,J}}(\omega)= F_{\mathcal{Q}(\{\mathrm{L}_{n} : n\in I\}\cup\{\mathrm{L}_{m}^{1} : m\in J\})}(\omega)$ and $\mathrm{L}_{1}$ is a subalgebra of $F_{\mathcal{V}_{I,J}}(\omega)$, by previous theorem it is enough to prove that $A\times \mathrm{L}_{1}\in \mathcal{Q}(F_{\mathcal{V}_{I,J}}(\omega))=\mathcal{Q}_{I,J}$ for every $A\in\mathcal{Q}(\{\mathrm{L}_{n} : n\in I\}\cup\{\mathrm{L}_{m}^{1} : m\in J\})$.
   Trivially $A\times \mathrm{L}_{1}$ is a bipartite member of $\mathcal{Q}(\{\mathrm{L}_{n} : n\in I\}\cup\{\mathrm{L}_{m}^{1} : m\in J\})=\mathcal{Q}_{I,J}^{1}$. Then, by Theorem \ref{BipQ}, $A\times \mathrm{L}_{1}\in \mathcal{Q}_{I,J}$ concluding the proof.  \hfill $\Box$

%
%

\section{Bases of admissible rules}

From the algebraization of $\L_{\infty}$ we obtain a 1 to 1 correspondence between quasivarieties of MV-algebras and finitary extensions of $\L_{\infty}$. In fact given a finitary extension $L$, $L$ is also algebraizable and its equivalent quasivariety semantics is its associated quasivariety. Viceversa if $K$ is a quasivariety of MV-algebras the logic $\models_{K} $ is its associated finitary extension, where $\models_{K} $  is defined as follows $\Gamma\models_{K}\varphi$ iff for every $A\in K$ and every evaluation $e: Prop(X)\to A$ if $e[\Gamma]=\{1\}$ then $e(\varphi)=1$.  Moreover  there is a translation from formulas to equations and a translation from equations to formulas that allows to obtain an quasiequational axiomatization of a quasivariety $K$ from the axiomatization of its associated finitary extension, and viceversa to get an axiomatization of a finitary extension $L$ from the quasiequational axiomatization of its equivalent quasivariety semantics.

Jerabek in \cite{Jer10b} gives an infinite  axiomatization for all $\L_{\infty}$-admissible rules ad moreover he proves that they are not finitely based. Our purpose is to obtain a base of all admissible rules for every proper axiomatic extension of $\L_{\infty}$. By Komori's classification of axiomatic extensions of $\L_{\infty}$ and Panti's  correspondence \cite{kom81, Pa99}, every axiomatic extension is given by a reduced pair $(I,J)$.  Given a reduced pair $(I,J)$ we denote by $L_{(I,J)}$ its associated axiomatic extension. Notice that $\mathcal{V}_{I,J}$ is the equivalent quasivariety semantics of $L_{I,J}$. Moreover since $\mathcal{V}_{I,J}=\mathcal{Q}(\{\mathrm{L}_{i} : i\in I\}\cup\{\mathrm{L}_{j}^{\omega} : j\in J\})$ (see \cite{GispThesis}), we get the following finite strong completeness theorem:

  \emph{$\varphi_{1},\ldots \varphi_{n} \vdash_{L_{I,J}} \varphi$ if and only if $\varphi_{1},\ldots \varphi_{n} \models_{\{\mathrm{L}_{i} : i\in I\}\cup\{\mathrm{L}_{j}^{\omega} : j\in J\}} \varphi$.}

\begin{lemma}\label{Rednm}

%

Let  $(I,J)$ be a reduced pair and $n= max\{max I, max J +1\}$. Then
$\lnot p^{m}\vdash_{L  _{I,J}} \lnot p^{n}$ for every $m>0$
\end{lemma}

\proof

By completeness $\lnot p^{m}\vdash_{L  _{I,J}} \lnot p^{n}$ is equivalent to  the following statement: \emph{For every  $A\in \{\mathrm{L}_{i} : i\in I\}\cup\{\mathrm{L}_{j}^{\omega} : j\in J\}$ and any $a\in A$, $a^{m}=0$ implies $a^{n}=0$}  which is valid because $Rad(A)=\{a\in A : a^{n}\neq 0\}$ for every $A\in \{\mathrm{L}_{i} : i\in I\}\cup\{\mathrm{L}_{j}^{\omega} : j\in J\}$. \hfill $\Box$

 \bigskip

In \cite{Jer10a} the author gives a basis of single conclusion passive rules for every extension of BL.

\begin{theorem}\label{PassAdmiss}{\rm \cite[Theorem 3.6]{Jer10a}}\\
If $L$ is an extension of BL then $CC^{1}=\{\lnot (p\lor\lnot p)^{n}/\bot : n>1\}$ is a basis of single-conclusion passive $L$-admissible rules.
\end{theorem}

\begin{theorem}
    Admissible rules for proper axiomatic extensions of $\L_{\infty}$ are finitely based.
 \end{theorem}
 \proof  By Theorem \ref{AllStComp}, every $L_{I,J}$-admissible rule is either derivable in $\models_{\mathcal{Q}_{I,J}^{1}}$ or it is a passive $L_{I,J}$-admissible rule. By Theorem \ref{PropQuasFinAx}, derivable rules in $\models_{\mathcal{Q}_{I,J}^{1}}$ are finitely axiomatizable. Moreover by Theorem \ref{PassAdmiss} and Lemma \ref{Rednm} passive $L_{I,J}$-admissible rules are axiomatizable by $\lnot (p\lor\lnot p)^{n}/\bot$ where $n= max\{max I, max J+1\}$. \hfill $\Box$

 \begin{theorem}\label{Ax}
   Let $(I,J)$ be a reduced pair, then a base of admissible rules for $\vdash_{I,J}$ is given by

\bi
   \item $\L1,\ \L2,\ \L3,\ \L4$ + M.P.
   \item $\alpha_{I,J}(\gamma)$.
   \item $\Delta(Q_{p}):= [(\lnot\varphi)^{p-1} \leftrightarrow \varphi]\lor[\psi\leftrightarrow\chi]\   / \ \psi\leftrightarrow\chi$\\ for every prime number $p\in Div(J)\smallsetminus Div(I)$

   \item $\Delta(U_{q}):=[(\lnot\varphi)^{q-1} \leftrightarrow \varphi]\lor[\psi\leftrightarrow\chi]\   / \ \alpha_{Iq,\emptyset}(\gamma) \lor (\psi\leftrightarrow\chi)$\\ for every prime number $q\in Div(I)$, where $I_{q}=\{n\in I : q|n\}$

   \item $CC_{n}^{1}:=\lnot(\varphi\lor \lnot \varphi)^{n}/\perp$\\ where $n= max\{max I, max J+1\}$
\ei

   \end{theorem}

\proof Following the proof of the previous theorem it is enough to put together the axiomatization of $\models_{\mathcal{Q}_{I,J}^{1}}$ plus $CC_{n}^{1}$ where $n= max\{max I, max J+1\}$. It follows from \cite[Theorem 4.5]{gispert} that
$\models_{\mathcal{Q}_{I,J}^{1}}$ is axiomatized by  

\bi
   \item $\L1,\ \L2,\ \L3,\ \L4$ + M.P.
   \item $\alpha_{I,J}(\gamma)$.
   \item $\Delta(Q_{n})$ for every  $n\in Div(J)\smallsetminus Div(I)$

   \item $\Delta(U_{m})$ for every  $m\in Div(I)$

   \item $Q_{n}:= (\lnot\varphi)^{p-1} \leftrightarrow \varphi\   / \ \psi$\\ for every  $n\in Div(J)\smallsetminus Div(I)$

   \item $U_{m}:=(\lnot\varphi)^{q-1} \leftrightarrow \varphi\   / \ \alpha_{Im,\emptyset}(\gamma) $ for every  $m\in Div(I)$
\ei

We can avoid $Q_{n}$ and $U_{m}$, since they  are passive $L_{I,J}$-admissible rules, therefore derivable from $CC_{n}^{1}$. It enough to take $\Delta(Q_{p})$ for every prime number $p\in Div(J)\smallsetminus Div(I)$, and $\Delta(U_{q})$ for every prime number $q\in Div(I)$ because $\Delta(Q_{n})$ is derivable from $\Delta(Q_{p})$ if $p|n$ and $\Delta(U_{m})$ is derivable from $\Delta(U_{q})$ if $q|m$.  \hfill $\Box$

\bigskip

\noindent \textbf{Acknowledgements:}  The author wishes to thank the anonymous referees for their suggestions that help to improve the presentation of this article. The   author is partially supported by grants  MTM2011-25747  of the  Spanish Ministry of Education and Science, including FEDER funds of  the European Union, the IRSES project MaToMUVI (PIRSES-GA-2009-247584) of the European Union and  2014SGR788  of Generalitat de Catalunya.

\bigskip

\noindent Joan Gispert. {\sc Facultat de Matem\`{a}tiques, Universitat de Barcelona, Gran Via 585, 08007 Barcelona Catalonia}\\ e-mail: \texttt{jgispertb@ub.edu}
\noindent

\end{document}